\documentstyle[11pt, amstex, amscd, %
epsfig, euscript, amssymb, amsfonts]{article}
\def\fg{{\frak{g}}}

\def\d{{\textsf{d}}}
\def\wdg{{{\mathfrak{w}}^{\d}_q(\mathfrak{g})}}
\def\cc{{\cal{C}}}
\def\lr#1{\langle #1\rangle}
\def\ol{\overline}
\def\a{\alpha}
\def\N{{\mathbb{N}}}

\def\otm{\otimes}
\def\Z{{\mathbb{Z}}}

\def\n{{\frak{n}}}
\def\p{{\verb"inj"}}
\def\aut{{\verb"Aut"}}
\def\id{\operatorname{\textsf{id}}}
\def\ch{{\mathcal{H}}}
\def\wgup{{{\mathfrak{w}}^{\ol{\d}}_q(\mathfrak{g})}}
\def\wgdd{{{\mathfrak{w}}^{\underline{\d}}_q(\mathfrak{g})}}

\def\db{\ol{\kappa}}
\def\dc{{\mathfrak{d}}}
\def\dbc{{\ol{\mathfrak{d}}}}

\def\bin#1#2{
\bigg[\begin{array}{c}
    #1\\
    #2\\
    \end{array}\bigg]
}
\newtheorem{theorem}{Theorem}[section]
\newtheorem{lemma}[theorem]{Lemma}
\newtheorem{proposition}[theorem]{Proposition}

\newtheorem{definition}[theorem]{Definition}
\newtheorem{example}[theorem]{Example}

\numberwithin{equation}{section}

\newenvironment{proof}{{\it\noindent Proof.}}{\hfill $\square$}

\date{}

\begin{document}

\pagestyle{myheadings} \markboth{\centerline{\sc s.
yang}}{ \centerline{\sc Weak quantum algebras}}

\title{\textsc{Weak Hopf algebras corresponding to Cartan matrices}}
\author{Shilin Yang\footnote{The  author is
partially supported by the National Science Foundation of China
(Grant. 10271014) and the Fund of Elitist Development
of Beijing City (Grant: 20042D0501518).}
\\
\small\sl College of Applied Sciences,
Beijing University of Technology\\
\small\sl Beijing 100022, P. R.  China \\
\small \verb" slyang@@bjut.edu.cn" }
\date{}
\maketitle

\begin{abstract}
We replace the group of
group-like elements of the quantized enveloping algebra $U_q(\fg)$
of a finite dimensional semisimple Lie algebra $\fg$
by some regular
monoid and get the weak Hopf algebra $\wdg$.
It is a new subclass of
weak Hopf algebras but not Hopf algebras.
Then we devote to constructing a basis of $\wdg$ and determine the group of weak Hopf
algebra automorphisms of $\wdg$ when $q$  is not a root of unity.
\end{abstract}

\section{Introduction}

 Recently, many mathematicians are interested in generalizations
 of Hopf algebras, of which importance has been recognized in both
 mathematics and physics. One way is to introduce the notion of a
 weak co-product, such that $\Delta(1)\ne 1\otm 1$,  which was
  motivated by the study of symmetries in low dimensional quantum
  field theory. This resulted in the definition of weak Hopf algebras,
 introduced by B$\ddot{o}$hm, Nill, and Szlach$\acute{a}$nyi (see for
example \cite{bns}). Since they are not bi-algebras, but almost bi-algebras,
there were also axioms required to define a weak antipode, differing
slightly from the usual notion of a Hopf algebra. The face algebras \cite{ha}
and generalized Kac algebras \cite{y} are examples of this class of
weak Hopf algebras.

It is possible to define a weak antipode
on a given bi-algebra. This was introduced by Li in \cite{fangli3}.
By definition, a bialgebra  $\ch=(H,\mu,\eta,\Delta,\varepsilon)$
over a field $k$
together with the identity map $\id$ in $\hom_k(H, H)$
is called a \emph{weak Hopf algebra} if there exists $T\in\operatorname*{hom}%
\nolimits_{k}(H,H)$ such that $\id\ast T\ast \id=\id$ and $T\ast \id\ast T=T$
where $\ast$ is the convolution product.
The map $T$ is called a weak antipode.
It is noted that
the notion of Hopf algebras, and left or right Hopf algebras
are included in this class of weak Hopf algebras
(see \cite{sweedler, nic/taf, gre/nic/taf}).
Another typical example is the weak quantum
algebras ${\frak{wsl}}_q(2)$ and $v{\frak{sl}}_q(2)$
constructed in \cite{ld}.  It is the generalization by replacing the set of
group-like elements of $U_q({\frak{sl}}_2)$ by the set of
some regular monoid, where  $U_q({\frak{sl}}_2)$
is the quantized enveloping algebra corresponding to 3-dimensional
semisimple Lie algebra.
The basis and some properties of ${\frak{wsl}}_q(2)$ (resp.
$v{\frak{sl}}_q(2)$) were
studied in \cite{ld}. Recently, Aizawa and Isaac \cite{ai} gave a description
of weak Hopf algebra ${\frak{wsl}}_q^d(n)$ in general, which
is corresponding to the known Hopf algebra $U_q({\frak{sl}}_n)$.

Our aim is to give more non-trivial examples
for weak Hopf algebras in the sense of Li.
Following
the idea \cite{ai, ld}, we would like to extend this construction
to the more general one $\wdg$ corresponding to  arbitrary  finite dimensional
semisimple
Lie algebra $\frak{g}$.

Thanks to  the definition of quantum group $U_q(\frak{g})$
defined by \cite{lu, jan},  we can also replace the group
$G(U_q(\fg))$ of
group-like elements by some regular
monoid and get the weak Hopf algebra $\wdg$,
which  is resulted from quantum group $U_{q}(\fg)$
\cite{ld, jan, lu}.
This successful construction provides us a new subclass of
weak Hopf algebras but not   Hopf algebras.  As does the classic
quantum
group $U_q(\fg)$, we will determine the basis and the group of weak Hopf
algebra automorphisms of $\wdg$.

To determine the basis of $\wdg$, we first show that  $\wdg$ can be written as a direct
sum of its two ideals and one of them is just isomorphic to the
classic quantum group $U_q(\fg)$. Then we apply the PBW Theorem
for $U_q(\fg)$  to describe  a basis of $\wdg$.
If $q$ is
not a root of unity, the group of Hopf algebra automorphisms
of $U_q(\fg)$ was determined in \cite{cm}.  The case when $q$ is a
root of unity, it was
considered in the recent work \cite{em}.  In the present paper,
we will determine the group of weak Hopf
algebra automorphisms of $\wdg$ under the condition that $q$ is not a root of unity.
The method is to apply the
result of [\ref{cm}, Corollary 4.3] and
some technical lemmas.

The paper is organized as follows.

In Sec. II, we give some
notations and the definition of weak Hopf algebra $\wdg$.
The ideal to construct the algebra $\wdg$ and some basic properties
 are described.  In Sec. III, we give the  comultiplication
of  $\wdg$ in order that it  is a weak hopf algebra but
not a Hopf alegbra.  The proof somewhat is basic and direct.
In Sec. IV, we describe the basis of $\wdg$ by the techonique
of  Lusztig's constructing PBW basis of
$U_q(\fg)$.
In the final section, we study and determine the group of Hopf automorphisms
of $U_q(\fg)$.

\section{Weak quantum algebras $\wdg$}

In this paper, we always assume that $k$ is a field of characteristic $0$.
Let $\fg$ be a finite dimensional semisimple Lie algebra. For the simplicity,
we can assume that $\fg$ is also simple. Then there is
a finite positive symmetrizable Cartan matrix $\cc=(a_{ij})_{n\times n}$
corresponding to it (see \cite{gtm9}).

Now we  let $R$ be the root system of $\fg$ and we fix a basis
$I=\{\a_1, \cdots, \a_n\}$ of $R$. Let $W$ be the Weyl group of $R$.
It is well known that
there is a unique $W$-invariant scalar
product $(\, ,\, )$ on the vector space generated by $R$ over the
reals such that $(\a, \a)=2$ for all short roots $\a$ in $R$.
Set for each $\a_i\in I, \ 1\leq i\leq n$
$$d_i=\frac{(\a_i, \a_i)}{2}.$$
It is noted that $(\a_i, \a_j)=d_i a_{ij}=d_j a_{ji}$.

Let $q\in k$ and  $q_i=q^{d_i}$, $1\leq i\leq n$.
It is assumed that $q_i\ne\pm1$, $0$ for all $i$. For an indeterminate $x$ and
an integer $m$,   let
 $$[m]_x=\frac{x^m-x^{-m}}{x-x^{-1}}, \  [m]!_x=[m]_x\cdots [1]_x, \
[0]_x!=1,$$
 and
$$\bigg[\begin{array}{c}
    {\small m}\\
    \small s\\
    \end{array}\bigg]_x=\frac{[m]!_x}{[s]!_x [m-s]!_x}.$$

One can review the definition of the quantized enveloping
algebra $U_q=U_q\left(\fg\right)$ by  referring to
\cite{kassel, jan, lu}. For the completeness, we describe it here
as follows. The algebra  $U_{q}(\fg)$ generated by  $4n$ generators
 $e_i$, $f_i$, $k_i$, $k_i^{-1} (1\leq i\leq n)$ with the relations%
\begin{eqnarray}
     &&k_ik_j=k_jk_i, \  k_i k_i^{-1}=k_i^{-1}k_i=1 \label{eqn1-1}\\
    &&k_ie_jk_i^{-1}=q_i^{a_{ij}}e_j,\
    k_if_jk_i^{-1}=q_i^{-a_{ij}}f_j; \label{eqn1-2}\\
    && e_if_j-f_je_i=\delta_{ij}\frac{k_i-k_i^{-1}}{q_i-q_i^{-1}};\label{eqn1-3}\\
    &&\sum_{s=0}^{1-a_{ij}} (-1)^s \bin{1-a_{ij}}{s}_{q_i}
    e_i^{1-a_{ij}-s} e_j e_i^s=0, \textrm{ if } i\ne j; \label{eqn1-4}\\
    &&\sum_{s=0}^{1-a_{ij}} (-1)^s
    \bin{1-a_{ij}}{s}_{q_i} f_i^{1-a_{ij}-s} f_j f_i^s
    =0, \textrm{ if } i\ne j. \label{eqn1-5}
\end{eqnarray}
To generalize the invertibility condition (\ref{eqn1-1}),
one way is to weak the invertibility to regularity, in which instead of
$\left\{  k_i, k_i^{-1}\right\}$ by a pair $\left\{  K_i, \overline{K}_i
\right\}$ for all $1\leq i\leq n$ subjecting to some relations. For example,
we can introduce a projector
$J$ such that
\begin{align}
&J=K_i\ol{K}_i=\ol{K}_i K_i, \nonumber\\
&J K_i =K_i,\;\;\;\;\;\;\;\overline{K}_{i}J=\overline{K}%
_{i},\nonumber\\
&J\overline{K}_{i}=\overline{K}_{i},\;\;\;\;\;\;\;K_{i}%
J=K_{i}.\nonumber%
\end{align}
for $1\leq i\leq n$.
To generalize other relations of definition, we need some terminologies
for simplicity. For example, if $E_i$  satisfies
\begin{eqnarray} \label{eqn-add-1}
K_j E_i=q_i^{a_{ij}} E_i K_j, \  E_i \ol{K}_j=q_i^{a_{ij}}\ol{K}_j E_i,
\ \forall j,
\end{eqnarray}
we say $E_i$ is type 1. Moreover, if $E_i$ satisfies
\begin{eqnarray} \label{eqn-add-2}
K_j E_i \ol{K}_j=q_i^{a_{ij}} E_i,
\ \forall j,
\end{eqnarray}
we say $E_i$  is type 2.
The same convention holds for $F_i$ by replacing $E_i$ with
$F_i$ and $a_{ij}$ with $-a_{ij}$ in the above relations.

We borrow some notations from the reference \cite{ai}, $2n$ simple generators $E_i$
and $F_i$ are listed by starting with the
$E_i$ followed by the $F_i$, where a $1$ is to indicate the use
of a type 1 generator and a $0$ is to the use of a type 2
generator. Then we write down a list of 0's and 1's
in the order corresponding to the generators determined by their type.
This gives a sequence $\d$ contained 0 and 1 in binary
representation of length $2n$. It is noted that
$\d$ contains all the information
on the relations with the generators $E_i$ and
$F_j$, all $K_j/\ol{K}_j$, and $J$.  We write $\d$
in terms of its binary expansion
$$\d=(\kappa_1,  \cdots, \kappa_n|\db_1, \cdots, \db_n),$$
where the bar separates the values representing the $E_i$ and $F_i$, and
where the $\kappa_i$  and $\bar{\kappa}_i$ have values of either 0 or 1. Accordingly, we can
say $E_i$ and $F_i$, $\a_i\in I$ are type $\d$ in an obvious sense.

Now we can write down this generalization explicitly as follows.
\begin{definition}
\label{def1} The algebra $\wdg$ is generated
by the $4 n+1$ variables $E_i$, $F_i$, $K_i$,
$\overline{K}_i (1\leq i\leq n)$ and $J$
with the relations: for all $1\leq i, j\leq n$,
\begin{align}
&J=K_i \ol{K}_i, \label{w0}\\
&K_i\overline{K}_j =\overline{K}_j K_i, K_i K_j=K_j K_i, \
\overline{K}_i\overline{K}_j=\overline{K}_j\overline{K}_i; \label{w1}\\
&J K_i=K_i,\;\;\;J\overline
{K}_i=\overline{K}_i, \label{w2}\\
&E_i, \ F_i \ \text{ \rm  are type  \d},\label{w3}\\
& E_i F_j-F_j E_i =\delta_{ij} \dfrac{K_i-\overline{K}_i}{q_i-q_i^{-1}};
\label{w5}\\
&\sum_{s=0}^{1-a_{ij}} (-1)^s \bin{1-a_{ij}}{s}_{q_i}
    E_i^{1-a_{ij}-s} E_j E_i^s=0, \textrm{ if } i\ne j; \label{w6}\\
&\sum_{s=0}^{1-a_{ij}} (-1)^s
    \bin{1-a_{ij}}{s}_{q_i} F_i^{1-a_{ij}-s} F_j F_i^s
    =0, \textrm{ if } i\ne j, \label{w7}%
\end{align}
The algebra $\wdg$ is said to be a $\d$-type \emph{weak quantum algebra} associated
to Lie algebra $\fg$.
\end{definition}

It is easy to see that there are $2^{2n}$ possible weak
quantum algebras $\wdg$ corresponding to the
sequence $\d$ in total.

It is easy to see that
(\ref{w1})--(\ref{w2})
generalize the relation (\ref{eqn1-1}), and
the relations (\ref{w3})--(\ref{w5}) generalize
the corresponding relations (\ref{eqn1-2})--(\ref{eqn1-3}).
The notations  $P_i$ ($1\leq i\leq n$) are defined by
$$P_i^k=\left\{%
\begin{array}{ll}
    K_i^k, &  \ k>0;\\
    J, &  \ k=0;\\
   {\ol{K}_i}^{-k}, & \ k<0. \\
\end{array}%
\right.
$$
It is easy to see that
$P_i^k$ satisfy the regularity conditions%
\begin{equation}
P_i^k P_i^{-k}P_i^k=P_t^k \label{jjj}%
\end{equation}
for all $k\in{\mathbb{Z}}$.

There are some  properties for $\wdg$ which are used later.

\begin{lemma}
\label{lem-center} The idempotent $J$ is in the centre of $\wdg$.
\end{lemma}

\begin{proof}
Indeed, for $K_j$ and $\ol{K}_j$,
it follows from (\ref{w1}) and (\ref{w2}).
For instance,
$$K_j J=K_j K_j \ol{K}_j=K_j\ol{K}_j K_j=J K_j.$$
For $E_j$, if it  is type 1, then we have
$$ J E_j=K_i\ol{K}_i E_j= q_i^{-a_{ij}} K_i E_j \ol{K}_i=
E_j J$$
by (\ref{eqn-add-1}); if $E_i$ is type 2,
$$JE_i=K_i\ol{K}_i E_i=q_i^{-2} K_i\ol{K}_i K_i E_i\ol{K}_i
=q_i^{-2} K_i E_i\ol{K}_i=q_i^{-2} K_i E_i\ol{K}_iK_i \ol{K}_i
=E_i J$$
by (\ref{eqn-add-2}). Hence, $JE_i=E_i J$ for all $i=1, 2, \cdots, n$.
The same argument for $F_i$ by (\ref{w3}).
\end{proof}

If $E_i$ is type 2, hence $K_j E_i\ol{K}_j=q_i^{a_{ij}} E_i$ for
all $1\leq j\leq n,$ then
$$K_jE_i=K_jJ E_i=K_jE_iJ=K_jE_i\ol{K}_j K_j=q_i^{a_{ij}} E_i K_j$$
and
$E_i \ol{K}_j=q_i^{a_{ij}}\ol{K}_j E_i$.
Similarly, if $F_i$ is type 2, we have
$$K_jF_i=q_i^{-a_{ij}} F_i K_j, \qquad F_i\ol{K}_j=q_i^{-a_{ij}} \ol{K}_j F_i.$$
Now, it is straightforward to check by induction that
$E_i$ (resp. $F_i$), $1\leq i\leq n$, is either type 1 or type 2,
the following relations hold in $\wdg$:
\begin{align}
&E_i^{m} K_j^{n}=q_i^{-mn a_{ij}}K_i^{n}E_i^{m},\;\;\;\;\;F_i^{m}%
K_j^{n}=q_i^{mn a_{ij}}K_j^{n}F_i^{m},\nonumber\\
&\label{ek1}\\
&E_i^{m}\ol{K}_j^{n} =q_i^{mn a_{ij}}\ol{K}_j^{n}E_i%
^{m},\;\;\;\;\;F_i^{m}\ol{K}_j^{n}=q_i^{-mn a_{ij}}\ol{K}_j^{n}%
F_i^{m}.\nonumber
\end{align}
In particular, we have
\begin{align}
&E_i^{m} K_i^{n}=q_i^{-2mn}K_i^{n}E_i^{m},\;\;\;\;\;F_i^{m}%
K_i^{n}=q_i^{2mn}K_i^{n}F_i^{m},\nonumber\\
&\label{ek2}\\
&E_i^{m}\ol{K}_i^{n} =q_i^{2mn}\ol{K}_i^{n}E_i%
^{m},\;\;\;\;\;F_i^{m}\ol{K}_i^{n}=q_i^{-2mn}\ol{K}_i^{n}%
F_i^{m}.\nonumber
\end{align}

\section{The weak Hopf algebra structure of $\wdg$ }

To make the $\d$-type weak quantum algebra $\wdg$  be a  weak Hopf algebra,
we define three maps
\begin{eqnarray*}
&&\Delta: \wdg\rightarrow \wdg\otimes
\wdg, \\
&&\varepsilon: \wdg\rightarrow k,\\
&&T: \wdg\rightarrow \wdg
\end{eqnarray*}
as follows:
\begin{align}
&\Delta(K_i)=K_i\otimes K_i,\\
&\Delta(\overline{K}_i)=\overline{K}_i\otimes\overline{K}_i,\\
&\Delta(J)=J\otimes J \label{dd2},\\
&\Delta (E_i) =
\left\{%
\begin{array}{ll}
   1\otimes E_i+E_i\otimes K_i, &  E_i \hbox{ is type 1;} \\
    J\otimes E_i+E_i \otimes K_i , & E_i \hbox{ is type 2,} \\
\end{array}%
\right. \\
&\Delta(F_i)=\left\{%
\begin{array}{ll}
F_i\otimes1+\overline{K}_i\otimes F_i, & F_i \hbox{ is type 1;} \\
  F_i\otimes J+ \ol{K}_i\otimes F_i, & F_i \hbox{ is type 2,} \\
\end{array}%
\right. \label{dd1}
\end{align}
$$
\varepsilon(E_i) =\varepsilon(F_i)=0,\;\varepsilon%
(K_i)=\varepsilon(\overline{K}_i)=1, \varepsilon(J)=1, \label{dd3}$$
while the map $T$ has the form
\begin{align}
T(1)&=1,\\
T(E_i)&=-E_i\overline{K}_i,\\
T_i(F_i)&=-K_iF_i\\
T(K_i)&=\overline{K}_i,\\
T(\overline{K}_i)&=K_i, T(J)=J.  %
\end{align}
Then we extend them to the whole $\wdg$.

In [\ref{ai}, Sec IV], the authors investigated the algebra
${\frak{wsl}}_q^{\textsf{d}}(n)$ and claimed that
${\frak{wsl}}_q^{\textsf{d}}(n)$ is a weak Hopf algebra, of which no
proof was given. In general, we yield that
\begin{theorem} \label{thm-m} {\rm ([\ref{ai}, Sec. IV])} \ $\wdg$ is  a non-commutative and
non-cocommutative weak Hopf algebra with the weak
antipode $T$, but \textit{not a Hopf algebra}.
\end{theorem}

In \cite{ai}, the authors gave a classification in some
sense of weak Hopf algebras corresponding to $U_q({\frak{sl}}_n)$. Similarly,
we can follow the idea \cite{ai} to describe the isomorphism classes
of weak Hopf algebras $\wdg$.
As a consequence, we have a lot of new non-trivial examples of
weak Hopf algebras for various sequences $\d$.

The theorem  follows from Lemma
\ref{lem2-1}  and Lemma \ref{lem2-2} below.

\begin{lemma} \label{lem2-1} $\wdg$ is a bialgebra  with
comultiplication $\Delta$ and counit $\varepsilon$.
\end{lemma}

\begin{proof}
It can be shown by direct calculation that the following relations hold.%
\begin{align*}
\Delta(K_i)\Delta(\overline{K}_j)&=\Delta(\overline{K}_j)\Delta(K_i),\\
\Delta(J)&=\Delta(K_i)\Delta(\ol{K}_i),\\
\Delta(J )\Delta(K_i)&=\Delta (K_i), \\
\Delta(J)\Delta(\ol{K}_i)&=\Delta(\overline{K}_i),\\
\varepsilon(K_i)\varepsilon(\overline{K}_j)&=\varepsilon
(\overline{K}_i)\varepsilon(K_j),\\ 
\varepsilon(J)\varepsilon(K_i)&=\varepsilon(K_i),\\ 
\varepsilon(J)\varepsilon (\overline{K}_i)&=\varepsilon(\overline{K}_i),\\
\varepsilon(\overline{K}_i)\varepsilon(E_i)&=q^{-2}%
\varepsilon(E_i)\varepsilon(\overline{K}_i),\\ 
\varepsilon(\overline{K}_i)\varepsilon(F_i)&=q^{2}%
\varepsilon(F_i)\varepsilon(\overline{K}_i),\\
\varepsilon(E_i)\varepsilon(F_j)-\varepsilon(F_j%
)\varepsilon(E_i)&=\delta_{ij}\dfrac{\varepsilon(K_i)-\varepsilon
(\overline{K}_i)}{q_i-q_i^{-1}};
\end{align*}
If $E_i$ is type 1, then
\begin{eqnarray*}
\Delta(K_j)\Delta(E_i)&=&(K_j\otimes K_j)(1\otimes E_i+E_i\otm K_i)\\
&=&K_j\otimes K_j E_i+K_j E_i\otimes K_j K_i\\
&=& q_i^{a_{ij}} K_j\otm E_i K_j+q_i^{a_{ij}} E_i K_j\otm K_i K_j\\
&=& q_i^{a_{ij}} \Delta(E_i)\Delta(K_j);
\end{eqnarray*}
if $E_i$ is type 2, then
\begin{eqnarray*}
\Delta(K_j)\Delta(E_i) \Delta(\ol{K}_j)
&=&(K_j\otm K_j)(J\otm E_i+E_i\otm K_i)(\ol{K}_j\otm \ol{K}_j)\\
&=&K_j \ol{K}_j\otm K_j E_i \ol{K}_j+K_j E_i \ol{K}_j \otm K_j K_i\ol{K}_j\\
&=&q_i^{a_{ij}} J\otm E_i+q_i^{a_{ij}} E_i\otm K_i\\
&=&q_i^{a_{ij}}\Delta(E_i).
\end{eqnarray*}
Therefore, $\Delta$ keeps the relation (\ref{w3}) for $E_i$'s.
The similar argument can show that $\Delta$ also keeps the
relation (\ref{w3}) for $F_i$'s.

Now we examine the identity
\begin{align*}
\Delta(E_i)\Delta(F_j)-\Delta(F_j)\Delta(E_i)  &
=\delta_{ij}\dfrac{\Delta(K_i)-\Delta(\overline{K}_i)}{q_i-q_i^{-1}}.
\end{align*}%
The following  cases should be considered:
 \begin{enumerate}
    \item $E_i$ is type 1 and $F_i$ is type 1,
    \item $E_i$ is type 1 and $F_i$ is type 2,
    \item $E_i$ is type 2 and $F_i$ is type 1,
    \item $E_i$ is type 2 and $F_i$ is type 2.
 \end{enumerate}
For the case 2, it is noted that
$$E_i \ol{K}_j=q_i^{a_{ij}} \ol{K}_j E_i=q^{d_i a_{ij}}\ol{K}_j E_i
=q^{(\alpha_i, \alpha_j)}\ol{K}_j E_i $$
and $$K_i F_j=K_iF_j J=K_i F_j \ol{K}_j K_i=q_j^{-a_{ji}} F_j K_i
=q^{-d_j a_{ji}} F_j K_i=q^{-(\alpha_i, \alpha_j)}F_j K_i.$$
The later identity holds since $K_iJ=K_i$ and $J$ is central in $\wdg$.
Then, it is easy to see that
\begin{eqnarray*}
\Delta(E_i)\Delta(F_j)-\Delta(F_j)\Delta(E_i)
&=&\ol{K}_j\otm(E_iF_j-F_j E_i)+(E_i F_j-F_j E_i)\otm K_i\\
&=&\delta_{ij} \ol{K}_j\otm \frac{K_i-\ol{K}_i}{q_i-q_i^{-1}}
+\delta_{ij} \frac{K_i-\ol{K}_i}{q_i-q_i^{-1}}\otm K_i\\
&=&\delta_{ij} \frac{K_i\otm K_i-\ol{K}_i\otm \ol{K}_i}{q_i-q_i^{-1}}
=\delta_{ij}\frac{\Delta(K_i)-\Delta(\ol{K}_i)}{q_i-q_i^{-1}}.
\end{eqnarray*}
We have shown that $\Delta$ keeps the relation (\ref{w5}) for the
case 2.  For the other cases  the proof is similar.
To see the map $\Delta$ keeps
the quantum Serre relations (\ref{w6}) and (\ref{w7}), we should consider
several cases according to the type of $\{E_i, E_j\}$ or $\{F_i, F_j\}$
($i\ne j$). In fact, for each case, the argument is more or less the same as
the case of $U_q(\fg)$ (see [\ref{jan}, p.67-68]).

Therefore,  $\Delta$ and
$\varepsilon$ can be extended to algebra morphisms from $\wdg$ to
$\wdg\otimes \wdg$ and from
$\wdg$ to $k$ respectively.

By the above relations it can be shown that
\begin{align*}
(\Delta\otimes\operatorname{id})\Delta(X)  &  =(\operatorname{id}%
\otimes\Delta)\Delta(X),\label{da}\\
(\varepsilon\otimes\operatorname{id})\Delta(X)  &  =(\operatorname{id}%
\otimes\varepsilon)\Delta(X)=X \label{ea}%
\end{align*}
for any $X=E_i, F_i, K_i$ or $\overline{K}_i$. Let $\mu$ and
$\eta$ be the product and the unit of $\wdg$
respectively, then $(\wdg, \mu, \eta, \Delta
,\varepsilon)$ becomes into a bialgebra.
\end{proof}

It is easy to see that
\begin{align*}
T(\overline{K}_i)T(K_j)  &  =T(K_j)T(\overline{K}_i),\label{d11}\\
T(J)T(K_i)  &  =T(K_i%
),\label{d12}\\
T(J)T(\overline{K}_i)  &  =T%
(\overline{K}_i),\label{d12a}\\
T(F_j)T(E_i)-T(E_i)T(F_j)  &  =\delta_{ij}\dfrac{T%
(K_i)-T(\overline{K}_i)}{q_i-q_i^{-1}}, \label{d15}%
\end{align*}
and $E_i$ is either type 1 or type 2, the map $T$ keeps
the anti-relation of (\ref{w3}). The argument for $F_i$ is
similar.  Moreover, $T$ also keeps the anti-relation for quantum
Serre relations.
For example, for $ 1\leq i, j\leq n$ with $i\ne j$, let
$r=1-a_{ij}$, we have
\begin{eqnarray*}
&&\sum_{s=0}^{r} (-1)^s
\bin{r}{s}_{q_i}
    T(E_i)^sT(E_j)T(E_i)^{r-s}\\
&&=\sum_{s=0}^{r} (-1)^s
\bin{r}{s}_{q_i} (-1)^{r+1} (E_i\ol{K}_i)^s (E_j\ol{K}_j)
(E_i\ol{K}_i)^{r-s}\\
&&=-q_j^2q_i^{r^2+r+ra_{ij}} \ol{K}_i^r \ol{K}_j
\sum_{s=0}^r (-1)^{r-s}
\bin{r}{s}_{q_i} E_i^s E_j E_i^{r-s}=0.
\end{eqnarray*}
Here we use the formula (\ref{ek1}) and (\ref{ek2}).
The argument for $F_i$ is similar. Therefore,
$T$ can be extended to an anti-algebra
morphism from $\wdg$ to $\wdg$ respectively.

Recall that the convolution product in the bialgebra $(\wdg, \mu, \eta,
\Delta, \varepsilon)$ is defined in the similar way to the
standard one(see e.g. \cite{kassel}) as
\begin{equation}
\left(  f\ast\, g\right)  \left(  X\right)  =\mu \left(  f\otimes
g\right)  \Delta(X) \label{ab}%
\end{equation}
for all $f, g\in \textrm{Hom}(\wdg, \wdg)$ and $X\in \wdg$.
It is noted that if $E_i$ is type 2, then
$$JE_i=K_i\ol{K}_i E_i=q_i^{-2} K_i\ol{K}_i K_i E_i\ol{K}_i
=q_i^{-2} K_i E_i\ol{K}_i=E_i.$$
The same argument shows that $F_i J=F_i J=F_i$ if $F_i$ is type 2.
Let $\id$ denote identity map in $\hom_k(\wdg, \wdg)$.
\begin{lemma} \label{lem2-2} Let $X$  be $E_i, \ F_i,\  K_i$, or $\overline{K}_i$,  then
\begin{align*}
(\id\ast \,T\ast\,\id)(X)  &
=\id(X),\label{it1}\\
(T\ast\,\id\ast\,T)(X)  &  =T(X).
\label{it2}%
\end{align*}
\end{lemma}

\begin{proof}
It is easy for $X=K_i, \overline{K}_i$.
We consider $X=E_i$, as an example. We set
$$\Delta_2=(\id\otm \Delta)\circ\Delta.$$
If $E_i$ is type 1, then
$$\Delta_2(E_i)=
1\otm 1\otm E_i+1\otm E_i\otm K_i+E_i\otm K_i\otm K_i.$$
It follows that
\begin{eqnarray*}
(\id\ast T\ast \id)(E_i)&=& T(1)E_i +T(E_i) K_i+E_i T(K_i) K_i\\
&=&E_i-E_i \ol{K}_i K_i+E_i \ol{K}_i K_i=\id(E_i),
\end{eqnarray*}
and
\begin{eqnarray*}
(T\ast \id\ast T)(E_i)&=& T(1)T(E_i) +T(1)E_i T(K_i) +T(E_i)K_i T(K_i)\\
&=&-E_i\ol{K}_i+E_i \ol{K}_i -E_i \ol{K}_i K_i \ol{K}_i\\
&=&-E_i \ol{K}_i=T(E_i).
\end{eqnarray*}
If $E_i$ is type 2, then
$$\Delta_2(E_i)=J\otm J\otm E_i+J\otm E_i\otm K_i+E_i\otm K_i\otm K_i.$$
It also deduce that
\begin{eqnarray*}
(\id\ast T\ast \id)(E_i)&=& J\; T(J)E_i + J T(E_i) K_i+E_i T(K_i) K_i\\
&=&J E_i-E_i \ol{K}_i K_i+E_i\ol{K}_i K_i=JE_i=\id(E_i)
\end{eqnarray*}
since $JE_i=E_i$,
and
\begin{eqnarray*}
(T\ast \id\ast T)(E_i)&=& T(J)J T(E_i) +T(J)E_i T(K_i) +T(E_i)K_i T(K_i)\\
&=&-E_i\ol{K}_i+E_i \ol{K}_i -E_i \ol{K}_i K_i \ol{K}_i\\
&=&-E_i \ol{K}_i=T(E_i).
\end{eqnarray*}
As for $F_i$, the argument is similar.
The proof of the lemma is finished.
\end{proof}

In order to conclude that the antipode axioms hold
on arbitrary elements, the following two facts are to be used.

   \quad (a)  The co-products of the generators are bilinear expressions
    of generators;

   \quad (b)  one of  $(\id\ast T)(X)$ and $T\ast \id(X)$ is a central
   element of $\wdg$ for $X$ being  the generators
   $K_i,  \ol{K}_i, E_i, F_i$.

The fact (a) is obvious. To see (b), we note the fact that
$(\id\ast T)(X)=\varepsilon(X)J$
for $X=K_i,  \ol{K}_i, E_i$ ($1\leq i\leq n)$ and $F_i$ of type 2.
Hence $(\id\ast T)(X)$ is in the center of $\wdg$.
However, if $X=F_i$
is type 1, the $(\id\ast T)(X)=(1-J)F_i$ may not be a central element, but
$(T\ast \id)(F_i)=\varepsilon(F_i)J$ is in the center.
Similarly,
$(T\ast\id)(X)$
for $X=K_i,  \ol{K}_i, F_i$ ($1\leq i\leq n)$ and $E_i$ is of type 2,
in the center of $\wdg$.
However, if $X=E_i$ is type 1, the $(T\ast\id)(X)=(1-J)E_i$ may not be a
central element, but
$(\id\ast T)(E_i)$ is in the center of $\wdg$.
This means that
(b)  holds.

It is noted that $E_i(1-J)F_j=F_j(1-J) E_i$ for all
$i, j\in \{1, 2, \cdots, n\}$ by the relation (\ref{w5}). Therefore,
if $E_i$ ( resp. $F_i$) is type 1, $(T\ast \id)(E_i)$ (resp. $(\id\ast T)(F_i)$)
commutates with all $F_j$( resp. all $E_j$), and $K_j, \ol{K}_j$
($1\leq j\leq n$).

We should show the claim that if
\begin{eqnarray*}
 &&(\id\ast T\ast \id)(x)=x, \ (T\ast \id\ast T)(x)=T(x);\\
 &&(\id\ast T\ast \id)(y)=y, \ (T\ast \id\ast T)(y)=T(y),
 \end{eqnarray*}
for all $x$ and $y$  being generators $E_i, F_i, K_i, \ol{K}_i$,
then
$$(\id\ast T\ast \id)(xy)=xy, \ (T\ast \id\ast T)(xy)=T(xy).$$
However, it is considerable direct by the above facts.
Now, that the antipode axioms hold on arbitrary elements is obvious
by induction.

If we assume that with the operations
$\mu, \eta, \Delta, \varepsilon$
the algebra $\wdg$ would possess an antipode $S$ so as to
become a Hopf algebra, then $S$ \;should satisfy $(S\ast\,\operatorname{id}%
)(J)=\eta\varepsilon(J)$, and it would follow that
$S(J)J=1$ and $J$ is invertible. It is impossible since $J(1-J)=0$.
This implies
that $\wdg$ is not a Hopf algebra with the above operators.
The proof of Theorem \ref{thm-m} is finished.   \hfill $\square$

It should be noted that if $\fg={\frak{sl}}_n$, the algebra $\wdg$
is just the mixtures ${\frak{wsl}}_q^d(n)$ in \cite{ai}. In particular,
 if $\fg={\frak{sl}}_2$,  $\wdg$ where $\d=(1|1)$ (resp.
 $\d=(0|0)$) coincides with
 ${\frak{wsl}}_q(2)$ (resp. ${v\frak{sl}}_q(2)$) given in \cite{ld}.

\section{The basis of $\wdg$}\label{sec4}

One can find the relationship between
$U_{q}(\fg)$
and the quantum algebra $\wdg$ as follows.
\begin{proposition}
\label{prop1} $\wdg/\lr{ J-1}\cong U_q(\fg)$.
\end{proposition}
\begin{proof}
 It is obvious by cancelling $K_i$.
\end{proof}

In fact, we can give a more  explicit relationship between
$\wdg$ and $U_q(\fg)$. For this purpose, we let
$w_q=\wdg J$ and $\ol{w}_q=\wdg (1-J)$.
We have the following decomposition.

\begin{proposition}\label{prop3-1}
\label{theor14}  As algebras $\wdg=w_q\oplus \ol{w}_q$. Moreover,
 $w_q\cong U_q(\fg)$ as Hopf algebras.
\end{proposition}
\begin{proof}
Noting that
$J$ is a central idempotent,
we see that $w_q$ and $\ol{w}_q$ are ideals of $\wdg$. It follows that
$$\wdg=w_q\oplus \ol{w}_q$$
as algebras.
Moreover, it is easy to see that
$w_q$ is generated by $E_i J, \ F_i J, \ K_i, \ \ol{K}_i$
and $J$ subject to the relations (\ref{w0})-(\ref{w2}) and
\begin{eqnarray}
&&K_i (E_jJ)=q_i^{a_{ij}} (E_j J) K_i,\;\;\;\overline{K}_i (E_j J)=
q_i^{-a_{ij}}(E_j J)\overline{K}_i,\label{eqn3-1}\\
&&K_i (F_j J) =q_i^{-a_{ij} }(F_j J) K_i,\;\;\;
\overline{K}_i (F_j J)=q_i^{a_{ij}}(F_j J)\overline{K}_i,\label{eqn3-2}\\
&&(E_iJ)(F_j J)-(F_j J)(E_i J) =\delta_{ij} \dfrac{K_i-\overline{K}_i}{q_i-q_i^{-1}};
\label{eqn3-3}\\
&&\sum_{s=0}^{1-a_{ij}} (-1)^s \bin{1-a_{ij}}{s}_{q_i}
    \left( E_i J\right )^{1-a_{ij}-s}(E_j J) \left(E_i J\right )^s=0,
    \textrm{ if } i\ne j; \label{eqn3-4}\\
&&\sum_{s=0}^{1-a_{ij}} (-1)^s
    \bin{1-a_{ij}}{s}_{q_i} \left( F_i J\right )^{1-a_{ij}-s}(F_j J)
    \left(F_i J\right )^s
    =0, \textrm{ if } i\ne j.\label{eqn3-5}%
\end{eqnarray}
Here $J$ can be viewed as the identity of $w_q$. At this point of view
$w_q$ becomes a Hopf algebra, in which the
co-multiplication $\Delta$ is
\begin{eqnarray*}
&&\Delta(E_i J)=J\otimes E_iJ + E_iJ\otimes K_i,\\
&&\Delta(F_i J)=F_i J\otimes J+\overline{K}_i\otimes F_i J,\\
&&\Delta(K_i)=K_i\otimes K_i,\
\Delta(\overline{K}_i)=\overline{K}_i \otimes\overline{K}_i.
\end{eqnarray*}
The counit $\varepsilon$ is
$$\varepsilon(E_i J)=\varepsilon(F_i J)=0,\ \varepsilon(K_i)
=\varepsilon(\overline{K}_i)=1.$$
The antipode $S$ is
$$S(E_i J)=-(E_i J)\overline{K}_i, \ S(F_i J)=-K_i(F_i J),\
S(K_i)=\overline{K}_i,\ S(\overline{K}_i)=K_i. $$
Let $\rho$ be the algebra morphism from $U_{q}(\fg)$ to $w_q$
defined by
$$\rho(e_i)=E_iJ,\ \rho(f_i)=F_i J,\ \rho(k_i)=K_i,\
\rho({k_i}^{-1})=\overline{K}_i.$$
It is straightforward to see that $\rho$ is a Hopf algebra
isomorphism.
\end{proof}

For the ideal $\ol{w}_q$ of $\wdg$, some conventions should
be noted.
Let
$$\d=(\kappa_1, \cdots, \kappa_n|\db_1, \cdots, \db_n)$$
be a binary sequence.
If $\kappa_i$ (resp. $\db_i$), $1\leq i\leq n$ is zero, and hence
that $E_i$ (resp. $F_i$) is type 2,
then $E_i(1-J)=0$ (resp. $F_i(1-J)=0$); if $\kappa_i$ (resp. $\db_i$),
$1\leq i\leq n$ is non-zero, and hence that $E_i$ (resp. $F_i$) is
type 1,  then $E_i(1-J)\ne 0$
(resp. $F_i(1-J)\ne 0$). Let
$$\dc=\{i | \kappa_i\ne 0\} \text{ and } \ \dbc=\{i | \db_i\ne 0\}$$
and
$$X_i=E_i(1-J) , \ Y_j=F_j(1-J),$$
where $i\in\dc$, $j\in \dbc$.
It is easy to see that
$\{X_i, Y_j | i\in\dc, \ j\in\dbc\}\cup\{1-J\}$ generate
the ideal $\ol{w}_q$ with
enjoying the first relation
\begin{eqnarray}
\label{relation-1}
&&X_iY_j=Y_j X_i,\ \textrm{ for all } i\in\dc, j\in\dbc
\end{eqnarray}
from the relation (\ref{w5}).

To see what other relations $X_i$ and $Y_i$ enjoy, we
consider the following two extreme cases

{\bf 1}. the case
$$\ol{\d}:=(\underset{n \textrm{copies}}{\underbrace{1, \cdots, 1}}|%
 \underset{n \textrm{copies}}{\underbrace{1, \cdots, 1}}).$$
In this case,  $\dc=\{1, \cdots, n\}$
and $\dbc=\{1, \cdots, n\}$. From the quantum Serre relations
(\ref{w6}) and (\ref{w7}),
we get that
\begin{eqnarray}
&&\sum_{s=0}^{1-a_{ij}} (-1)^s \bin{1-a_{ij}}{s}_{q_i}
    \left( X_i\right )^{1-a_{ij}-s} X_j \left(X_i\right )^s=0,
    \textrm{ if } i\ne j, \label{eqn3-6}\\
&&\sum_{s=0}^{1-a_{ij}} (-1)^s
    \bin{1-a_{ij}}{s}_{q_i} \left( Y_i\right )^{1-a_{ij}-s} Y_j
    \left(Y_i\right )^s
    =0, \textrm{ if } i\ne j, \label{eqn3-7}
\end{eqnarray}
and other relations corresponding to (\ref{w0})-(\ref{w3}) would be vanished
automatically. This means that
the ideal $\ol{w}_q$ can be understood as an algebra generated
by $X_i, \ Y_i$, $1\leq i\leq n$. with an identity $1-J$ subject  to the relations
(\ref{relation-1})-(\ref{eqn3-7}).

{\bf 2}. the case
 $$\underline{\d}:=(\underset{n \textrm{copies}}{\underbrace{0, \cdots, 0}}|%
 \underset{n \textrm{copies}}{\underbrace{0, \cdots, 0}}).$$
In this case, $\dc$ and $\dbc$ are empty and
$\ol{w}_q=k(1-J)$.

In general, three cases should be considered.

\begin{enumerate}
    \item $\dc\ne \emptyset$ and $\dbc=\emptyset$;
    \item $\dc=\emptyset$ and $\dbc\ne\emptyset$;
    \item $\dc\ne \emptyset$ and $\dbc\ne\emptyset$.
\end{enumerate}
In the first case, as an algebra  $\ol{w}_q$ is generated
by $X_i, \ i\in\dc$ with an identity $1-J$ subject  to the relations
(\ref{eqn3-6}) with $i, j\in\dc$. In the second case, as an algebra
$\ol{w}_q$ is generated by $Y_i, \ i\in\dbc$ with an identity $1-J$
subject  to the relations (\ref{eqn3-7}) with $i, j\in\dbc$.
In the third case, $\ol{w}_q$ can be viewed as an algebra generated
by $X_i, \ Y_j$, $i\in\dc, \ j\in\dbc$ with an identity $1-J$
subject  to the relations (\ref{relation-1})-(\ref{eqn3-7}).

To consider the PBW basis of $\wdg$, we need some knowledge of
braid groups.  We define a simple reflection $s_i$ by
$$s_i(\a_j)=\a_j-a_{ij}\a_i$$
for all $i$ and $j$.
Let $W$ be the Weyl group of $R$; it is the subgroup of
$\textrm{GL}({\mathbb{Z}}^n)$ generated by the refelctions $s_i (1\leq i\leq n)$.
Let $\ell(w)$ be the usual length function on $W$ with respect
to the generators $\{s_1, \cdots, s_n\}$. Let
$R^+$ be the set of positive roots of $R$ with respect to the set
of simple roots $\Pi$ and $\ell_0=|R^+|$. For each pair
 $1\leq i, \, j\leq n$ with $i\ne j$, we let $r=-a_{ij}.$

As is known in \cite{jan}, for each $1\leq i\leq n$, there is a unique
automorphism  $T_i: U_q(\fg)\to U_q(\fg)$ such that
\begin{eqnarray*}
&&T_i(e_i)=-f_i k_i,\ T_i(f_i)=-k_i^{-1} e_i,\\
&&T_i(k_\mu)=k_{s_i(\mu)},\\
&&T_i(e_j)=\sum_{k=0}^r (-1)^i q_i^{-k} e_i^{(r-i)} e_j e_i^{(i)},\\
&&T_i(f_j)=\sum_{k=0}^r (-1)^i q_i^k f_i^{(i)} f_j f_i^{(r-i)}.
\end{eqnarray*}
They are called Lusztig's symmetries.
It is well known that $\{T_i | 1\leq i\leq n\}$ satisfies the braid
relations, that is
\begin{eqnarray*}
&& T_i T_j T_i=T_j T_i T_j, \quad \textrm{ if } s_i s_j \textrm{ of order } 2,\\
&& T_i T_j T_i T_j =T_j T_i T_j T_i, \quad \textrm{ if } s_i s_j
\textrm{ of order } 4,\\
&&T_i T_j T_i T_j T_i T_j=T_j  T_i T_j T_i T_j T_i,
\quad \textrm{ if } s_i s_j \textrm{ of order } 6.
\end{eqnarray*}
Therefore, the above facts allow us to define for each $w\in W$
an automorphism
$T_w$ of $U_q(\fg)$ as follows. For $w=1$ set $T_1=1$ (the identity). For
$w\ne 1$ choose a reduced expression
$w=s_{i_1}\cdots s_{i_m}$ and set
$$T_w=T_{i_1}\cdots T_{i_m}.$$
It is independent of the reduced expression.
Let $U_q^+$ (resp. $U_q^-$ and $U_q^0$) be the subalgebra of $U_q(\fg)$
generated by $e_i$ (resp. $f_i$ and $k_i, \, k_i^{-1}$), $1\leq i\leq n$.
Let $w_0$ be the
longest element in $W$ and let $w_0=s_{i_1}\cdots s_{i_t}$ be
a reduced expression.
Let $\mathbb{N}$ be the set of non-negative numbers.
According to this order we denote
$a=(a_t, \cdots, a_1)\in {\mathbb{N}}^{\ell_0}$
and
\begin{eqnarray}\label{eqn3-9}
e^a=T_{i_1}\cdots T_{i_{t-1}}(e_{i_t}^{a_t})\cdots
 T_{i_1} T_{i_2} (e_{i_3}^{a_3}) T_{i_1} (e_{i_2}^{a_2}) e_{i_1}^{a_1}.
\end{eqnarray}

The following theorem is well known.

\begin{theorem}{\rm (cf. [\ref{jan}, Theorem 8.24])}\label{lem3-2}
 The elements
$e^a$ {\rm (resp. $f^a$)} for $a\in {\mathbb{N}}^{\ell_0}$, are linearly
independent  and a basis of $U_q^+$ {\rm (resp. $U_q^-$).}
\end{theorem}

We note that the multiplication map
\begin{eqnarray}\label{eqn3-10}
U_q^-\otimes U_q^0\otimes U_q^+\to U_q(\fg), \quad
u_1\otimes u_2\otimes u_3\to u_1 u_2 u_3
\end{eqnarray}
is an isomorphism of vector spaces.

For $s=(s_1, s_2, \cdots, s_n)\in\Z^n$,
we define
$$P^s=P_1^{s_1} P_2^{s_2}\cdots P_n^{s_n}.$$

First, let us examine two examples.
\begin{example}\label{exam1}
The set
$$\left\{ F^b  P^s
E^a J\;\left
|\right. a, b\in{\mathbb{N}}^{\ell_0}, s\in \Z^n \right\}
\bigcup \left\{ F^b\,E^a(1-J)
\left |\right. a, b\in{\mathbb{N}}^{\ell_0} \right\}$$
forms a basis of $\wgup$.
\end{example}
\begin{proof}
Let $w_q^0$ be the subalgebra generated by $\{K_i, \ol{K}_i \left|\right.
1\leq i\leq n\}.$ It is easy
to see that $P^s$ ($s\in\Z$) forms a basis of $w_q^0$.

Let $w_q^+$ (resp. $w_q^-$) denote the subalgebra generated by
$E_iJ$(resp. $F_iJ$), $1\leq i\leq n$.

We replace $e_{i_k}$ where $k=1, \cdots, t$ in the
right hand side of (\ref{eqn3-9})
by $E_{i_k}J$ (resp. $E_{i_k}(1-J)$),  and the corresponding
left hand side
by $\left (E J\right)^a$ (resp. $\left(E(1-J)\right)^a$).
By Theorem \ref{lem3-2} the set $\left\{\right. \left(EJ\right)^a \left
|\right. a\in \N^{\ell_0}\left\}\right.$ (resp.
$\left\{\right. \left(FJ\right)^b \left
|\right. b\in\N^{\ell_0}\left\}\right.$)
forms a basis of $w_q^+$ (resp. $w_q^-$). It is easy to see
that $$\left(FJ\right)^b P^s
\left(EJ\right)^a =F^b P^s E^a J.$$
It follows from (\ref{eqn3-10}) that
$$\left\{ F^a  P^s
E^b J\;\left
|\right. a, b\in{\mathbb{N}}^{\ell_0}, s\in \Z^n \right\}$$
forms a basis of $w_q$.

Similarly, $\left\{ F^aE^b(1-J)
| a, b\in{\mathbb{N}}^{\ell_0} \right\}$ forms a basis of $\ol{w}_q$.
The proof is completed by Proposition \ref{prop3-1}.
\end{proof}

In a similar way we can get that
\begin{example}\label{exam2}
The set
$$\left\{ F^b  P^s
E^a J\;\left
|\right. a, b\in{\mathbb{N}}^{\ell_0}, s\in \Z^n \right\}
\bigcup \left\{1-J \right\}$$
forms a basis of $\wgdd$.
\end{example}

In general, let both $i$ and $j$ be in $\dc$ or  in $\dbc$,
we say that $i\thicksim j$ if there exist some sequence
$i=\gamma_1, \cdots,
\gamma_p=j$ in $\dc$,  where $\gamma_p$ is in $\dc$ or $\dbc$, $p$
is some positive integer,
such that
$(\a_{\gamma_i}, \a_{\gamma_{i+1}})\ne 0$  for all $i=1, \cdots, p-1$.
This is an equivalent relation.
Let $\frak{s}$ and $\ol{\frak{s}}$
be the set of equivalent classes on $\dc$
and $\dbc$ respectively.
If $\frak{e}$ is an element in $\frak{s}$ or $\ol{\frak{s}}$,
it is obvious that
$C_{\frak{e}}=(a_{ij})_{i, j\in \frak{e}}$ is also a symmetrizable
Cartan matrix.
If $\frak{e}_1\ne \frak{e}_2$ in $\frak{s}$ or $\ol{\frak{s}}$,
and $i\in \frak{e}_1$ and $j\in \frak{e}_2$, then $E_iE_j=E_j E_i$ and
$F_i F_j=F_j F_i$,
and hence $X_i X_j=X_j X_i$ and $Y_i Y_j=Y_j Y_i$ respectively.
Let $W_{\frak{e}}$ be the Weyl group corresponding to
the equivalent class $\frak{e}$ on $\dc$ or $\dbc$.
Let $w_0^{\frak{e}}=s_{i_1} \cdots s_{i_{t_{\frak{e}}}}$ be the
longest element in $W_{\frak{e}}$. Let
$\ell_{\frak{e}}=|w_0^{\frak{e}}|=t_{\frak{e}}$.
We denote $a_{\frak{e}}=(a_{t_{\frak{e}}}, \cdots, a_1)$
according to this order
and
\begin{eqnarray}\label{eqn3-12}
X_{\frak{e}}^{a_{\frak{e}}}&=&T_{i_1}\cdots T_{i_{t_{\frak{e}}-1}}
(E_{i_{t_{\frak{e}}}}^{a_{t_{\frak{e}}}})\cdots
 T_{i_1} T_{i_2} (E_{i_3}^{a_3}) T_{i_1} (E_{i_2}^{a_2}) E_{i_1}^{a_1}(1-J),\\
 Y_{\frak{e}}^{a_{\frak{e}}}&=&T_{i_1}\cdots T_{i_{t_{\frak{e}}-1}}
(F_{i_{t_{\frak{e}}}}^{a_{t_{\frak{e}}}})\cdots
 T_{i_1} T_{i_2} (F_{i_3}^{a_3}) T_{i_1} (F_{i_2}^{a_2}) F_{i_1}^{a_1}(1-J).
\end{eqnarray}
It is noted that
$$X_i Y_j=Y_j X_i$$
for all $i\in\dc, j\in\dbc$, and
$$\dc=\bigcup_{\frak{e}\in\frak{s} } {\frak{e}},\quad
\dbc=\bigcup_{\frak{e}\in\ol{\frak{s}}} {\frak{e}},
$$
one sees that
$$\left\{\prod_{\frak{e}\in \frak{s}} X_{\frak{e}}^{a_{\frak{e}}}
\prod_{\ol{\frak{e}}\in\ol{\frak{s}}} Y_{\ol{\frak{e}}}^{b_{\ol{\frak{e}}}}
| a_{\frak{e}}\in \N^{\ell_{\frak{e}}},
b_{\ol{\frak{e}}}\in \N^{\ell_{\ol{\frak{e}}}} \right \}$$
forms a basis of $\ol{w}_q$.  By Proposition \ref{prop3-1} and the discussion
above, it follows that
\begin{theorem}\label{thm3-4}
The notations are kept as above. Then the set
$$\left\{ F^b  P^s
E^a J\;\left
|\right. a, b\in{\mathbb{N}}^{\ell_0}, s\in \Z^n \right\}
\bigcup \left\{\prod_{\frak{e}\in \frak{s}} X_{\frak{e}}^{a_{\frak{e}}}
\prod_{\ol{\frak{e}}\in\ol{\frak{s}}} Y_{\ol{\frak{e}}}^{b_{\ol{\frak{e}}}}
\;|\;  a_{\frak{e}}\in \N^{\ell_{\frak{e}}},
b_{\ol{\frak{e}}}\in \N^{\ell_{\ol{\frak{e}}}}
\right \}$$
forms a basis of $\wdg$.
\end{theorem}

It is mentioned that the set
$$\left\{ E^a  P^s
F^a J\;\left
|\right. a, b\in{\mathbb{N}}^{\ell_0}, s\in \Z^n \right\}
\bigcup \left\{\prod_{\frak{e}\in \frak{s}} X_{\frak{e}}^{a_{\frak{e}}}
\prod_{\ol{\frak{e}}\in\ol{\frak{s}}} Y_{\ol{\frak{e}}}^{b_{\ol{\frak{e}}}}
\;|\;  a_{\frak{e}}\in \N^{\ell_{\frak{e}}},
b_{\ol{\frak{e}}}\in \N^{\ell_{\ol{\frak{e}}}}
\right \}$$
also forms a basis of $\wdg$.

Let us recall some basic facts used below.
Let $C$ be a coalgebra. If the set $\{C_n\}_{n\ge 0}$ of subspaces
 of $C$ satisfies
\begin{enumerate}
    \item $C_n\subseteq C_{n+1} \textrm{ and  }  C=\cup_{n\ge 0} C_n$
    \item $\Delta(C_n)\subseteq \sum_{i=0}^n C_i\otm C_{n-i},$
\end{enumerate}
then the set  $\{C_n\}_{n\ge 0}$  is said to be
a coalgebra filtration of $C$.

The following lemma is well known.
\begin{lemma} \label{lem3-4}{\rm ([\ref{mon}, Lemma 5.5.1])}
Let $H$ be a bialgebra which contains subspaces $A_0\subset A_1$
such that
\begin{enumerate}
    \item  $A_0$ is a (unital) subalgebra of $H$ and $A_1$ is
    a left, and a right $A_0$-module;
    \item  $A_1$ generates $H$ as an algebra, and $1\in A_0$;
    \item  $\Delta A_0\subseteq A_0\otm A_0$ and $\Delta(A_1)\subseteq A_1\otm A_0
    +A_0\otm A_1$.
\end{enumerate}
Then, if we set $A_n=(A_1)^n$ for all $n\ge 1$, $\{A_n\}$ is
a coalgebra filtration of $H$ and $A_0\supseteq H_0$, where
$H_0$ is the coradical of $H$.
 \end{lemma}

The element $x\in \wdg$ is said to be a group-like element if
$\Delta(x)=x\otimes x.$
The set $G=G(\wdg)$ of
all group-like elements of $\wdg$ can be  determined.

A semigroup $S$ is called regular, if for every $x\in S$, there
exists a $y\in S$ such that $xyx=x$ and $yxy=y$ and a monoid is
a semigroup with identity.
\begin{proposition}\label{prop3-5}
The set of all group-like elements is $G=\{P^s\left |\right.
s\in\Z^n\}\bigcup \{1\}$, which forms a regular monoid under the
multiplication of $\wdg$.
\end{proposition}

\begin{proof}
Let $G$ be the setup as above. Let $A_0=kG$ and
$$A_1=A_0\left(\sum_i kE_i+kF_i+A_0\right )A_0.$$
 It is easy to see that  $A_0\subset A_1$ satisfies the hypotheses
of Lemma \ref{lem3-4} and $A_0\subseteq H_0$. Hence $H_0=A_0$.
\end{proof}

It is mentioned that $\wdg$ is a pointed bialgebra
with the coradical $kG$ by Proposition \ref{prop3-5}.

\section{The automorphism group of {\rm $\wdg$}}

If $(A, m, \mu, \Delta, \varepsilon, T)$ is a (weak) Hopf algebra, then
a (weak) Hopf algebra automorphism $\varphi: A\to A$  is an invertible algebra
homomorphism satisfying
\begin{eqnarray*}
&&(\varphi\otm \varphi)\circ \Delta=\Delta\circ \varphi,\\
&& \varepsilon=\varepsilon\circ \varphi,\\
&&\varphi\circ T=T\circ \varphi.
\end{eqnarray*}
The group of Hopf algebra automorphisms of $U_q(\fg)$ was
determined by several authors. See for example \cite{cm},
\cite{em}.  Inspired by these considerations, we would like to
determine the group of weak Hopf algebra automorphisms of $\wdg$
where $q$  is not a root of unity.

let $N=(k^*)^n$, and for $a=(a_1, \cdots, a_n)\in N$, we define a
map $\phi_a: U_q(\fg)\to U_q(\fg)$ by
$$\phi_a(k_i)=k_i, \ \phi_a(e_i)=a_i e_i, \ \phi_a(f_i)=a_i^{-1} f_i.$$
It is straightfoward to check that $\phi_a$ is a Hopf algebra
automorphism of $U_q(\fg)$. It is called $N$ the group of diagonal
automorphisms of $U_q(\fg)$.

 Recall that the Dynkin diagram of $\fg$ is the weight graph $\Gamma$
 with vertices $\n=\{1, 2, \cdots, n\}$ such that vertices $i$ and $j$
 are connected by $a_{ij}a_{ji}$ edges, and vertex $i$ carries weight
 $d_i$. Let $\sigma$ be a automorphsim of Dynkin diagram $\Gamma$, that
 is
 $\sigma$ is a bijection of $\n$ and
 $$d_i a_{ij}=d_{\sigma(i)} a_{\sigma(i)\sigma(j)}$$
 for all $1\leq i, j\leq n$. If this is the case, there is an automorphism
 of Hopf algebra $U_q(\fg)$, also denote by $\sigma$, given by
 $$\sigma(k_i)=k_{\sigma(i)}, \ \sigma(e_i)=e_{\sigma(i)}, \
 \sigma(f_i)=f_{\sigma(i)}.$$

 We denote by $H$ the group of automorphisms of the Dynkin diagram.
 Also, $H$ acts on $N$ by the rules $\sigma\cdot a=(a_{\sigma(1)}, \cdots
 , a_{\sigma(n)})$ and we have
 $\phi_a\sigma=\sigma \phi_{\sigma\cdot a}$. We will base on the following
 theorem to investigate the group of automorphisms of $\wdg$.

\begin{theorem} \label{lem4-2}
{\rm ([\ref{cm}, Corollary 4.3; \ref{tw}, Theorem 2.1])} The
group of Hopf algebra automorphism of $U_q(\fg)$ is the semidirect
product $N\rtimes H$ of the group of diagonal automorphism $N$ by
the group of diagram automorphisms $H$.
\end{theorem}

Moreover, we also need some basic lemmas.
\begin{lemma}\label{lem4-3}
If $x, \ y\in \ol{w}_q$, $a_x,\ b_y\ne 0$ and
\begin{eqnarray*}
&&\Delta(x)=1\otm x+x\otm K_i+ a_x (1-J)\otm E_i J,\\
&&\Delta(y)=\ol{K}_i\otm y+y\otm 1+ b_y F_i J\otm (1-J).
\end{eqnarray*}
then $x=a_x E_i (1-J)$ and $y=b_y F_i(1-J)$.
\end{lemma}
\begin{proof} The notations $\dc$ and $\dbc$ are as in Section \ref{sec4}.
Let $W_{\dc}$ be the Weyl group  corresponding to $\dc$ and
$w_{\dc}=s_{i_1}\cdots s_{i_{t_\dc}}$ be
the longest element in $W_{\dc}$ and $\ell_{\dc}=t_{\dc}$.
Similarly, we have the Weyl group $W_{\dbc}, \;
w_{\dbc}=s_{j_1}\cdots s_{j_{t_{\dbc-1}}}, \; \ell_{\dbc}$
for the support set $\dbc$ in an obvious sense.

Let $\wdg^+$ be the sub-bialgebra of $\wdg$ generated by $E_i$, $1\leq i\leq n$
 and the set $G$ of group-likes of $\wdg$.
Let $\wdg^-$ be the sub-bialgebra of $\wdg$ generated
by $F_i$, $1\leq i\leq n$ and $G$.  We define a
$\N[I]$-algebra gradation on $\wdg^+$ (resp. $\wdg^-$)
 such that $E_j^s$ (resp. $F_j^s$) are homogeneous of degree
 $s\a_j\in\N[I]$ for $s\in\N$, $1\leq j\leq n$. We also set
 $\deg K_j=\deg \ol{K}_{j}=\deg J=0$ for all $j$.
 According to this gradation, $\wdg^+$ (resp. $\wdg^-$) is also a graded
 coalgebra. It is obvious that $\wdg^+(1-J)\subset \ol{w}_q$
 has a basis
 $\{E^a(1-J)|a\in\N^{\ell_\dc}\}$. Similarly,
 $\wdg^-(1-J)\subset \ol{w}_q$  has a basis
 $\{F^b(1-J)|b\in\N^{\ell_\dbc}\}$.

It is easy to see the elements $E^a(1-J)$, $F^b(1-J)$ have a
gradation $|a|$ and $|b|$, where
\begin{eqnarray*}
|a|&=&a_1\a_1+a_2 s_{i_1}(\a_{i_2})+\cdots+a_{t_{\dc}} s_{i_1}s_{i_2}
\cdots s_{i_{t_{\dc}-1}}(\a_{i_{t_{\dc}}})\\
|b|&=&b_1\a_1+b_2 s_{j_1}(\a_{j_2})+\cdots+b_{t_{\dc}} s_{j_1}s_{j_2}
\cdots s_{j_{t_{\dbc}-1}}(\a_{i_{j_{\dbc}}}).
\end{eqnarray*}

Let $$x=\sum_{a\in\N^{\ell_{\dc}}, b\in\N^{\ell_{\dbc}}}x(a, b) E^a F^b (1-J)\in \ol{w}_q,$$
where $\{E^a F^b (1-J)\}$ are linear independent.  We have
\begin{eqnarray*}
\Delta(x)&=&\sum_{a\in\N^{\ell_{\dc}},\; b\in\N^{\ell_{\dbc}}\atop{|s|+|t|=|a|\atop{|x|+|y|=|b|}}}
x(a, b)
a(s, t) b(x, y) E^s F^x\ol{K}_{|y|}\otm E^t K_{|s|} F^y\\
&\quad&-\sum_{a\in\N^{\ell_{\dc}}, b\in\N^{\ell_{\dbc}}
\atop{|s|+|t|=|a|\atop{|x|+|y|=|b|}}}x(a, b)
a(s, t) b(x, y) E^s F^x\ol{K}_{|y|}J\otm E^t K_{|s|} F^y J.
\end{eqnarray*}
Hence,
\begin{eqnarray*}
\Delta(x)
&=&\sum_{a\in\N^{\ell_{\dc}}, b\in\N^{\ell_{\dbc}}\atop{|x|+|y|=|b|\atop{|a'|=|a|, |y|\ne 0}}}
x(a, b)a(0, a')b(x, y) F^x\ol{K}_{|y|}\otm E^{a'} F^y(1-J)\\
&\quad&+\sum_{a\in\N^{\ell_{\dc}}, b\in\N^{\ell_{\dbc}}
\atop{|s|+|t|=|a|\atop{|b'|=|b|, |s|\ne 0}}}
x(a, b)a(s, t)b(b', 0)
E^s F^{b'}(1-J)\otm  E^t K_{|s|}\\
&\quad &+\sum_{a\in\N^{\ell_{\dc}}, b\in\N^{\ell_{\dbc}}
\atop{|a'|=|a|\atop |b'|=|b|}}x(a, b)
a(0, b') b(b', 0) F^{b'}\otm E^{a'}\\
&\quad&-\sum_{a\in\N^{\ell_{\dc}}, b\in\N^{\ell_{\dbc}}\atop{|a'|=|a|\atop |b'|=|b|}}
x(a, b)a(0, b')b(b', 0) F^{b'}J\otm E^{a'}J.
\end{eqnarray*}
On the other hand, by the assumption  we have
\begin{eqnarray*}
\Delta(x)&=&1\otm \sum_{a\in\N^{\ell_{\dc}}, b\in\N^{\ell_{\dbc}}}x(a, b) E^a F^b(1-J)
+\sum_{a\in \N^{\ell_{\dc}}, b\in\N^{\ell_{\dbc}}}x(a, b) E^a F^b(1-J)\otm K_i
\\
\qquad &+& a_x (1-J)\otm E_iJ.
\end{eqnarray*}
Comparing the above identities, we conclude that
all $b=0$. Now we can rewrite $x$ as
$$x=\sum_{a\in\N^{\ell_{\dc}}}x(a) E^a (1-J)\in \ol{w}_q,$$ and
\begin{eqnarray*}
&&\sum_{a\in\N^{\ell_{\dc}}\atop{|s|+|t|=|a|\atop{|s|\ne 0}}}
x(a)a(s, t) E^s(1-J)\otm  E^t K_{|s|}+\sum_{a\in\N^{\ell_{\dc}}\atop
{|a'|=|a|}}x(a) 1\otm E^{a'}-\sum_{a\in\N^{\ell_{\dc}}\atop
{|a'|=|a|}} x(a) J\otm E^{a'} J\\
&&=\sum_{a\in\N^{\ell_{\dc}}\atop{|s|+|t|=|a|\atop{|s|\ne 0}}}
x(a)a(s, t) E^s(1-J)\otm  E^t K_{|s|}+
1\otm \sum_{a\in\N^{\ell_{\dc}}\atop{|a'|=|a|}}x(a) E^{a'}(1-J)\\
&&+\sum_{a\in\N^{\ell_{\dc}}} x(a)(1-J)\otm E^{a'}J\\
&&=1\otm \sum_{a\in\N^{\ell_{\dc}}}x(a) E^a(1-J)
+\sum_{a\in\N^{\ell_{\dc}}}x(a) E^a(1-J)\otm K_i+
a_x(1-J)\otm E_iJ.
\end{eqnarray*}
Also, comparing the above identity, we conclude
that $|a|=|a'|=i$,  and hence $x=x(a)E_i(1-J)$.
It  follows that $x(a)=a_x$. The argument  for $F_i$ is similar.
\end{proof}
\begin{lemma}\label{lem4-4}
Let $\varphi$ be a Hopf algebra automorphism of $w_q$ (the identity is $J$).
Then there exists a unique way to extend $\varphi$ to $\wdg$
such that $\varphi$ is an  automorphism of weak Hopf algebra $\wdg$.
\end{lemma}
\begin{proof}
Let $\varphi$ be the automorphism of Hopf algebra $w_q$. By
Lemma \ref{lem4-2},
the map $\varphi$ is
$$K_i\to K_{\sigma(i)}, \ \ol{K}_i\to \ol{K}_{\sigma(i)}, \
E_iJ\to a_{\sigma(i)} E_{\sigma(i)}J, \
F_i J\to a_{\sigma(i)}^{-1} F_{\sigma(i)}J$$
for some automorphism $\sigma$ of Dynkin diagram.
Assume that $\varphi$ can be extended to $\wdg$ such that $\varphi$
is an automorphism of $\wdg$ as weak Hopf algebras. We have to find a suitable
images of $\varphi(E_i(1-J))$ and $\varphi(F_i(1-J))$.  For example,
if $E_i$  is type 2, we do nothing  since
$E_i(1-J)=0$.  Assume that $E_i$ is type 1 and
$\varphi(E_i(1-J))=x$, then $x\ne 0$, $x\in \ol{w}_q$, and
$$\varphi(E_i)=\varphi(E_iJ+E_i(1-J))=a_{\sigma(i)} E_{\sigma(i)}J+x.$$
Since $\varphi$ is a coalgebra homomorphism, we have
\begin{eqnarray*}
&&a_{\sigma(i)} E_{\sigma(i)}J \otm K_{\sigma(i)}+a_{\sigma(i)}
J\otm E_{\sigma(i)}J+\Delta(x)\\
&=&
(a_{\sigma(i)} E_{\sigma(i)}J+x) \otm K_{\sigma(i)}+
1\otm (a_{\sigma(i)} E_{\sigma(i)}J+x).
\end{eqnarray*}
It follows that
$$\Delta(x)=x\otm K_{\sigma(i)}+1\otm x+
a_{\sigma(i)}(1-J)\otm E_{\sigma(i)}J.$$
By Lemma \ref{lem4-3}, $x=a_{\sigma(i)}E_{\sigma(i)}(1-J).$
Similarly, we can get that $\varphi(F_i(1-J))=
a_{\sigma(i)}^{-1} F_{\sigma(i)}(1-J)$
if  $F_i$ is type 1, and nothing is done if $F_i$ is of type 2.
Of course, $\varphi(1)=1$. The lemma is proved.
\end{proof}


For $a=(a_1, \cdots, a_n)\in N$, we define a map
$\phi_a: \wdg\to \wdg$ by
$$\phi_a(K_i)=K_i, \phi_a(\ol{K}_i)=\ol{K}_i,\ \phi_a(E_i)=a_i E_i, \ \phi_a(F_i)=a_i^{-1} F_i.$$
It is straightfoward to check that $\phi_a$ is a weak Hopf algebra
automorphism of $\wdg$.
If $\sigma\in H$, then there is an automorphism
 of weak algebra $\wdg$, also denote by $\sigma$, given by
 $$\sigma(K_i)=K_{\sigma(i)}, \ \sigma(\ol{K}_i)=\ol{K}_{\sigma(i)}, \
 \sigma(E_i)=E_{\sigma(i)}, \
 \sigma(F_i)=F_{\sigma(i)}.$$
Recall that there is an action of $H$ on $N$ by the rules
$\sigma\cdot a=(a_{\sigma(1)}, \cdots , a_{\sigma(n)})$ and
 $\phi_a\sigma=\sigma \phi_{\sigma\cdot a}$.

Let $\aut(\wdg)$ be the group of automorphisms of  weak Hopf algebra
$\wdg$. The group $\aut(\wdg)$ can be determined by the following theorem.
\begin{theorem}\ $\aut(\wdg)= N\rtimes H$.
\end{theorem}

\begin{proof}
Theorem A in \cite{cm} is  the key to determine the group
of automorphisms of Hopf algebra. But we don't know whether
it is true or not for bialgebras (see [\ref{cm1}]), we can not
directly apply Theorem A to get the result.

Let $\varphi\in\aut (\wdg)$, since $\varphi$ sends group-likes to group-likes,
we have
 $\varphi(J)=P^s$ for some
$s\in\Z^n$ by Proposition \ref{prop3-5}. If $s\ne 0$, since $J^2=J$
we have $P^{2s}=P^s$, hence $P^s=J$ it follows that $P^sJ=J$.
It is a contradiction for Theorem \ref{thm3-4}. Therefore, $s=0$ and
$\varphi(J)=J$.

According to Proposition \ref{prop3-1},
$\wdg=w_q\oplus \ol{w}_q$ and \ $w_q\cong U_{q}(\fg)$ as Hopf algebras,
where the notations $w_q$ and $\ol{w}_q$ as before.

Let $\p_q: w_q\to \wdg$ be the inclusion  defined by
\begin{eqnarray*}
J\to J, \ E_i J\to E_i J, \ F_i J\to F_i J, \ K_i\to K_i, \
\overline{K_i} \to \overline{K}_i,
\end{eqnarray*}
and then extend it by linearity. It is easy to see that $\p_q$ is
a weak Hopf algebra injection by Proposition \ref{prop3-1}. Let
$\varphi\in\aut(\wdg)$, we see that $w_q=\textrm{im}(\varphi\circ
\p_q)$ since $\varphi(J)=J$. This implies that $\varphi\circ \p_q$
is an automorphism of Hopf algebra $w_q$ and $\varphi\circ\p_q\in
N\rtimes H$ by Theorem \ref{lem4-2}. That is, the map
$\varphi|_{w_q}$ is
$$K_i\to K_{\sigma(i)}, \ \ol{K}_i\to \ol{K}_{\sigma(i)}, \
E_iJ\to a_{\sigma(i)} E_{\sigma(i)}J, \
F_i J\to a_{\sigma(i)}^{-1} F_{\sigma(i)}J.$$
This implies that
$\varphi$ is
$$K_i\to K_{\sigma(i)}, \ \ol{K}_i\to \ol{K}_{\sigma(i)}, \
E_i\to a_{\sigma(i)} E_{\sigma(i)}, \
F_i \to a_{\sigma(i)}^{-1} F_{\sigma(i)}$$
 by Lemma \ref{lem4-3}.
Hence $\varphi\in N\rtimes H$ and
$\aut(\wdg)\subseteq N\rtimes H$. On the other hand, it
is obvious that $N\rtimes H\subseteq \aut(\wdg)$.
The proof is completed.
\end{proof}

The more interesting  problem is to determine the algebraic group of
$\wdg$. It is mentioned that
\cite{AD} contributed  to determine
$\textrm{Aut}_{\rm Alg}(U_q^+({\mathfrak g}))$.
A basic idea to approach the group $\textrm{Aut}_{\rm
Alg}\,(U_q^+({\mathfrak g}))$ is to study its actions on natural
sets.

\section*{Acknowledgment}

\quad The author is grateful to the referee for his/her helpful
comments.


\end{document}